\documentstyle[12pt]{article}

\def\CC{{\rm\kern.24em\vrule width.02em height1.4ex depth-.05ex\kern-.26em C}} 
\def\QQ{{\rm\kern.24em\vrule width.02em height1.4ex depth-.05ex\kern-.26em Q}}
\def\PP{{\rm\kern.24em\vrule width.02em height1.4ex depth-.05ex\kern-.26em P}}
\def\Rr{{\rm I\kern-.2em R}}
\def\ZZ{{\rm\kern.26em\vrule width.02em height0.5ex depth0ex\kern.04em\vrule width.02em height1.47ex depth-1ex\kern-.34em Z}}
\def\BB{{\rm\kern.24em\vrule width.02em height1.4ex depth-.05ex\kern-.26em B}}
\def\RR{\hspace{.065in}\rm{\vrule width.02em height1.55ex depth-.07ex\kern-.3165em R}}

\newcommand{\bea}{\begin{eqnarray*}}
\newcommand{\eea}{\end{eqnarray*}}

\def\CQFD{\hfill \vrule width 7pt height 7pt depth 1pt}

\newtheorem{theorem}{Theorem}

\newtheorem{main}{Main Theorem}
\newtheorem{lemma}{Lemma}
\newtheorem{proposition}{Proposition}

\newtheorem{problem}{Problem}
\newtheorem{problem 1}{Problem 1}
\newtheorem{problem 2}{Problem 2}
\newtheorem{problem 3}{Problem 3}
\newtheorem{definition}{Definition}
\newtheorem{remark}{Remark}

\begin{document}





\centerline{The Julia set of Henon maps}



\medskip

\centerline{{\scshape John Erik Forn\ae ss}\footnote{July 2, 2005. The author is supported by an NSF grant}}
{\footnotesize
 \centerline{Mathematics Department, University of Michigan}
  \centerline{530 Church Street, 2074 East Hall}
     \centerline{Ann Arbor, MI 48109-1043, USA}
   }



 \medskip

\begin{abstract}
 In this paper we investigate the support of the unique measure of maximal entropy
 of  complex H\'enon maps, $J^*$. The main question is whether this set is the same as the analogue of the Julia set
 $J$.
 \end{abstract}

\section{Introduction}

Let $H:\CC^2 \rightarrow \CC^2$ be an automorphism of the form
$H(z,w)=(P(z)+aw,bz)$ where $P:\CC \rightarrow \CC$ is a polynomial
of degree $d$ at least $2$ and $a,b$ are nonzero complex constants. More generally, we define a
(complex) H\'enon map $H:\CC^2 \rightarrow \CC^2$
as $H=H_1 \circ \cdots \circ H_n$ where each $H_i=(P_i(z)+a_iw,b_iz)$ is of the
above form and $n \geq 1.$

\medskip

Associated to each H\'enon map there is a natural invariant measure $\mu$ with
compact support $J^*$. There is also a natural notion of Julia set, $J.$
We recall the precise definitions in Section 2. The following is one of the basic questions
in the theory of complex dynamics.

\begin{problem}
Is $J=J^*$?
\end{problem}

It was proved in \cite{BS1991} that if $H$ is uniformly hyperbolic when restricted to $J$ then
$J=J^*$. This leaves open the following interesting special case of the Problem:

\begin{problem} If $H$ is uniformly hyperbolic on $J^*$, is $J=J^*?$
\end{problem}

For motivation, we recall that this question arose naturally in the author's investigation of sustainable complex H\'enon maps, see \cite{F2002} which also contains several references to the dynamics of complex H\'enon maps
and also includes an exposition of background material for this article. The author showed that
complex H\'enon maps which are uniformly hyperbolic on $J$ are sustainable and that sustainable
maps are uniformly hyperbolic on $J^*.$ The last step in the characterization of sustainable
complex H\'enon maps was to show the equality of $J$ and $J^*$, however, this used
sustainability. We are left with the equivalent problem:

\begin{problem}
If $H$ is uniformly hyperbolic on $J^*,$ is $H$ sustainable?
\end{problem}

We also note in passing that it is an interesting open question whether there is a similar characterization
of sustainable real H\'enon maps.

\medskip

\begin{main}
Let $H$ be a complex H\'enon map which is hyperbolic on $J^*.$ If
$H$ is not volume preserving, then $J=J^*.$
\end{main}

In Section 2 we introduce notation and review background results. We prove the Main
Theorem in Section 3.
The author would like to thank Eric Bedford for valuable comments.

\section{Notation and background results}

We recall some standard notation for H\'enon maps which can be found in many sources, see
for example \cite{F2002}. Let $H$ be a complex H\'enon map. We denote by
$H^n$ the $n$- fold composition $H \circ \cdots \circ H$ 
 for any positive
integer $n \geq 1$. If $n<0$ we write $H^n:=(H^{-1})^{|n|}$ 
where $H^{-1}$ denotes the inverse
map. Also $H^0 ={\mbox{Id}}.$ We define the sets of  bounded orbits and their boundaries.

\bea
K^+ & :=  & \{(z,w)\in \CC^2; \{H^n(z,w)\}_{n \geq 1} \; {\mbox{is a bounded sequence}}\}\\ 
K^- & :=  & \{(z,w)\in \CC^2; \{H^n(z,w)\}_{n \leq -1} \; {\mbox{is a bounded sequence}}\}\\
K & :=  & K^+ \cap K^-\\
J^+ & :=  & \partial K^+\\
J^- & :=  & \partial K^-\\
J & :=  & J^+ \cap J^-\\
\eea

 We let $d$ denote the degree of the highest order term in the polynomial mapping
 $H, d=\Pi_{i=1} {\mbox{deg}}(P_i).$ We define the escape functions $G^\pm:\CC^2
 \rightarrow \RR$ by
 
 \bea
 G^+ (z,w)& :=  \lim_{n \rightarrow \infty} \frac{\log^+ \|H^n(z,w)\|}{d^n}\\
 G^- (z,w)& :=  \lim_{n \rightarrow -\infty} \frac{\log^+ \|H^n(z,w)\|}{d^n}\\
 \eea
 
 We remark that the sets $K^\pm$ are closed and the functions $G^\pm$ are continuous
 and plurisubharmonic. Moreover $G^\pm$ vanishes on $K^\pm$ and is strictly positive
 and pluriharmonic on $\CC^2 \setminus K^\pm.$
 
 \medskip
 
 We define $\mu^\pm:= dd^c G^\pm$ and $\mu:=\mu^+ \wedge \mu^-.$ The $(1,1)$ currents
 $\mu^\pm$ are supported on $J^\pm$ and hence $\mu$ is a positive measure supported inside $J$. 
 We set
 $J^*$ equal to the support of $\mu.$ So we get immediately that $J^* \subset J.$
 
\medskip

We recall the notion of uniform hyperbolicity. Let $F:M \rightarrow M$ be a
holomorphic automorphism of a complex manifold $M$ of dimension $m$ with a 
metric $d.$ Suppose that the compact set $S,$ $S \subset M$ is  completely invariant, i.e.
$F(S)=S.$ 

\medskip

We say that $F$ is uniformly hyperbolic on $S$ if there exists a continuous 
splitting $E^s_x \oplus E^u_x=T_x, x\in S$ of complex subspaces of the
tangentspaces. Moreover $F'(E^s_x)=E^s_{F(x)}$ and $F'(E^u_x)=E^u_{F(x)}.$
Also there exist constants $0<c, \lambda<1$ so that
$\|(F^n)'(v)\|\leq \frac{\lambda^n}{c}\|v\|$ for every 
$v\in E^s_x, x\in S, n\geq 1$ and $\|(F^n)'(v)\| \geq \frac{c}{\lambda^n}
\|v\|$ for every $v\in E^u_x, x\in S, n \geq 1$.

\medskip

For $p\in S,$ let $W^s_p:= \{x\in M;
d(F^n(x),F^n(p))\rightarrow 0, n \rightarrow \infty\}$
and $W^u_p:=\{x\in M; d(F^n(x),F^n(p))\rightarrow 0, n \rightarrow 
-\infty\}$. Then $W^s_p$ and $W^u_p$ are immersed complex manifolds
of the same dimension as $E^s_p$ and $E^u_p$ respectively. Similarly,
we define 
$W^s_S:=\{x\in M; d(F^n(x),S) \rightarrow 0, n\rightarrow \infty\}$ 
and 
 $W^u_S:=\{x\in M; d(F^n(x),S) \rightarrow 0, n\rightarrow -\infty\}$,
the stable and unstable sets of $S.$

\begin{lemma} \cite{BS1991B} Let $\tilde{M}$ be a one dimensional complex manifold
in $\CC^2.$ Let $M \subset \subset \tilde{M}$ be smoothly bounded. If
$\mu^+_{| \tilde{M}}(\partial M)=0$, then 
$$\frac{1}{d^n}H^n_*([M]) \rightarrow c \mu^-, c=\mu^+_{|\tilde{M}}(M).$$
\end{lemma}

Using this Lemma and the arguments of Lemma 9.1 and Theorem 9.6 in \cite{BLS1993}, it follows that
for some large $n,$ the manifold $H^n(M)$ must intersect the stable manifold
of any given periodic saddle point in $J^*$ transversely in any given neighborhood
of $J^*.$ Hence the same is true for $M.$
[There might also be some tangential intersections.]
So we obtain:

\begin{proposition}
Suppose that $H$ is a complex H\'enon map. Let $\tilde{M}$ be a one dimensional complex submanifold of $\CC^2$ and let
$M\subset \subset \tilde{M}$ such that $\mu^+_{|\tilde{M}}(M)>0.$ Let $U\supset J^*$
be some neighborhood.  Then for any saddle point $p$ there are arbitrarily large integers
so that the manifold $H^n(M)$
intersects $W^s(p)$ transversally in $U.$
\end{proposition}

\begin{proposition}
Suppose that $H$ is a complex H\'enon map. Let $p\in J \setminus J^*.$ Then there are arbitrarily
small discs $M \subset \tilde{M}$ through $p$ on which $\mu^+_{|\tilde{M}}(M)>0.$
For every transverse intersection $x$ near $p$ between $M$ and a stable manifold of
a saddle point, the function $G^-$ vanishes identically in a neighborhood
of the intersection in the stable manifold provided $G^-(x)=0,$
so in particular if $G^-$ vanishes identically on $M,$ which happens
if $M$ is an unstable manifold.
\end{proposition}

{\bf Proof:}
If not, then the local piece of stable manifold must contain
a transverse intersection $y$ with the unstable manifold of some other saddle
point. But then the $\alpha$ and $\omega$ limit sets of $y$ are periodic saddles.
This implies that $y$ is a limit of saddle points, so must be in $J^*,$ a contradiction.\\

\CQFD

It is convenient to state the same results for the inverse map.

\begin{proposition}
Suppose that $H$ is a complex H\'enon map. Let $\tilde{M}$ be a one dimensional complex submanifold of $\CC^2$ and let
$M\subset \subset \tilde{M}$ such that $\mu^-_{|\tilde{M}}(M)>0.$ Let $U\supset J^*$
be some neighborhood.  Then for any saddle point $p$ there are arbitrarily large integers $n$ 
so that the manifold $H^{-n}(M)$
intersects $W^u(p)$ transversally in $U.$ 
\end{proposition}

\begin{proposition}
Suppose that $H$ is a complex H\'enon map. Let $p\in J \setminus J^*.$ Then there are arbitrarily
small discs $M \subset \tilde{M}$ through $p$ on which $\mu^-_{|\tilde{M}}(M)>0.$
For every transverse intersection $y$ near $p$ between $M$ and an unstable manifold of
a saddle point, the function $G^+$ vanishes identically in a neighborhood
of the intersection in the unstable manifold, provided that $G^+(y)=0$, so in 
particular if $G^+$ vanishes identically on $M$, which happens if $M$ is a
stable manifold.
\end{proposition}

\begin{theorem}
Let $H$ be any complex H\'enon map. Then $J^*$ is equal to the closure of the
collection of periodic saddle orbits.
\end{theorem}

This is in \cite{BLS1993}.

\section{Proof of the Main Theorem}

 We start by recalling a few standard lemmas which can for example be found in
\cite{F2002}.

 \begin{lemma}
 If $H$ is uniformly hyperbolic on $J^*,$
 then $W^s_{J^*}= \cup_{p\in J^*} W^s_p.$
 \end{lemma}
 
 This is \cite{F2002}, Lemma 4.15.

 \begin{lemma}
 Suppose $H$ is uniformly hyperbolic on $J^*$ and that $p\in W^s_{J^*} \cap W^u_{J^*}.$
 Then $p\in J^*.$
 \end{lemma}

This is \cite{F2002}, Lemma 4.17.

 \begin{lemma}
 If $H$ is uniformly hyperbolic on $J^*$, then $J^*$ has local product structure.
 \end{lemma}
 
 This is \cite{F2002}, Lemma 4.14.

\medskip

We investigate some consequences of
assuming that $J \setminus J^*$ is nonempty. We give first some notation.
\noindent Let $p \in J^*$  and let $\Phi^s_p: \CC \rightarrow  \CC^2$
be a parametrization of the stable manifold $W^s_p $ of $p,$
$\Phi^s_p(0)=p$.  This exists a.e. d$\mu$ on $J^*$ (Oseledec regular points)
and if $H$ is uniformly hyperbolic on $J^*,$
this exists for all $p\in J^*.$ Also, let $J^*_{p,s}:= (\Phi^s_p)^{-1}(J^*)
$.\\

For the unstable manifold we use the analogous notation, $q, \Psi^u_q, J^*_{q,u}.$
 
\medskip

We next define a notion describing how $J\setminus J^*$ might be attached to $J^*.$ 
Namely we give a name to points whose stable manifolds or unstable manifolds
belong (partly) to $J\setminus J^*.$ We will say that a point in $J^*$ is stably
exposed (to $J \setminus J^*$) if its stable manifold enters into $J\setminus J^*$ and similarly for
unstably exposed. More precisely: 

\begin{definition}
We say that $p$ is stably exposed if $0\in\CC$ is a boundary point of a simply connected open set 
$V_{s,p} \subset \CC$ such that
$\partial V_{s,p} \subset J^*_{p,s}$
and $\Phi^s_p(V_{s,p})
\subset J \setminus J^*.$
 If we instead have 
$\Phi^s_p(V_{s,p})\subset K\setminus J^*,$ we
say that $p$ is weakly stably exposed.
\end{definition}

We define unstably exposed similarly, using the notation
$q, W_{u,q}.$

\begin{remark}
When the Jacobian of $H$ has modulus $1$, $W^s_p,W^u_q$ cannot intersect the interiors of
$K^\pm$ so weakly (un)stable is equivalent to (un)stable. If the modulus is $\neq 1$, either $K^+$
or $K^-$ has empty interior. 
\end{remark}

\begin{theorem}
Suppose that $H$ is a complex H\'enon map which is uniformly hyperbolic on $J^*.$ If $J \setminus J^*$
is nonempty, then there are points $p,q\in J^*$ so that $p$ is weakly stably
exposed and $q$ is weakly unstably exposed.
\end{theorem}

\begin{definition}
Suppose that $H$ is hyperbolic on $J^*$. Let $p\in J^*$ and let $\Phi^s_p:\CC\rightarrow
W^s_p, \Psi^u_p:\CC \rightarrow W^u_p$ parametrize the stable and unstable
manifolds of $p$, $\Phi^s_p(0)=\Psi^u_p(0)=p.$
 Suppose that $0\in \CC$ is a boundary point of simply connected components
  $U\subset \CC\setminus J^*_{p,s}$ and
  $V \subset \CC \setminus J^*_{p,u}.$
  In this case, because of local product structure,
 $J^*$ is locally contained on one side of each of the laminated local
  hypersurfaces which are the local unstable set of $\Phi^s_p(\partial U)$ and the
  local stable set of $\Psi^u_p(\partial V)$ respectively. We say in this case
  that $p$ is a distinguished boundary point of $J^*.$
  \end{definition}
  
  \begin{remark}
  Heuristically one can think of the local stable set as
  $\{(z,w)\in \CC^2; |w|=1\}$ and the local unstable set
  as $\{(z,w)\in \CC^2; |z|=1\}$ and that $J^*$ contains
  the distinguished boundary of the bidisc and $J^*$ is contained in the
  bidisc.
  \end{remark}

 \begin{theorem}
 Suppose that $H$ is a complex H\'enon map which is uniformly hyperbolic on $J^*.$ If $J \setminus J^*$
is nonempty, then distinguished boundary points are dense in $J^*.$
\end{theorem}

 \begin{proposition}
 Suppose that $H$ is uniformly hyperbolic on $J^*.$ 
 Let $p\in J^*$ and let $\Phi^s_p:\CC\rightarrow \CC^2$ parametrize
 the stable manifold of $p,$ $\Phi^s_p(0)=p.$ 
 Fix a nonempty connected component $V$ of $\CC\setminus J^*_{p,s}.$
 The function $G^+ \circ \Phi^s_p \equiv 0$ on $\CC.$ On the other hand
 $(G^-\circ \Phi^s_p)_{|V}$ is either identically zero or strictly positive.
 \end{proposition}
 
 {\bf Proof:}
 The first part is obvious. Suppose there is a point $z\in V$ so that $G^-\circ \Phi^s_p(z)=0$ but
 that $G^-\circ \Phi^s_p$ is not identically zero on $V.$ Choosing another $z$ if necessary we can
 assume that there are  discs $z\in \Delta \subset \subset \tilde{\Delta}\subset \subset V$ so that
 $G^-\circ \Phi^s_p(z)=0$ but that $G^-\circ \Phi^s_p$ is not vanishing identically on $\Delta.$
 Then the nonnegative function $G^-\circ \Phi^s_p$ cannot be harmonic on $\Delta.$ Let
 $M:=\Phi^s_p(\Delta)$. Then $\mu^-_{|M}$ has mass on $M$. Let $U$ be a neighborhood of
 $J^*$ and let $w \in J^*$ be a saddle point. By Proposition 3
 there exists an integer $n>>1$ so that $H^{-n}(M)$ intersects
 the unstable manifold of $w$ transversally in $U.$  Let $x$ be such an intersection
 point. Then $x \in W^s_p \cap W^u_w.$ By Lemma 3 this implies that $x\in J^*.$ Hence
 $H^n(x)\in J^*.$ This contradicts that $(\Phi^s_p)^{-1}(H^n(x))\in \Delta \subset \CC
 \setminus (\Phi^s_p)^{-1}(J^*).$\\
 
 \CQFD

Without the hypothesis of uniform hyperbolicity, the proof still works if we assume that
$p$ is a periodic saddle point. Hence:

\begin{proposition}
Let $H$ be a complex H\'enon map. Let $p$ be a periodic saddle point.
Suppose $V$ is a connected component of $\CC\setminus J^*_{p,s}.$
Then $G^-_{|V}$ either vanishes identically or is strictly positive.
\end{proposition}

{\bf Proof:}
Suppose $x \in V, \Delta \subset V, x\in \Delta$ and
$G^-(x)=0$ while $G^-$ does not vanish identically on $\Delta.$
Then there is a $y\in \Delta$ so that $y$ is a transverse intersection
with an unstable manifold of a periodic saddle point.
But then $y\in J^*$, a contradiction.\\

\CQFD

The following result is obvious.

\begin{lemma} Suppose that $H$ is uniformly hyperbolic
on $J^*.$ Let $p,q\in J^*.$ Suppose that $q\in W^u_p.$ Then $W^u_q=W^u_p.$
Likewise if $q\in W^s_p, $ then $W^s_q=W^s_p.$
\end{lemma}

We can now prove Theorem 2.

{\bf Proof:}
 Suppose that for every
$p\in J^*, G^-$ is not vanishing on any open subset of $W^s(p).$ Then there are
open subsets $J^* \subset U \subset \subset V$ so that $H(V)
\supset \supset U$, $J \setminus \overline{V} \neq \emptyset$ and $G^->c>0$
on $W^s_{loc}(J^*)\cap (V \setminus U).$ Pick $x\in J \setminus \overline{V}$. Then
$x\in J^+$ and hence, if $y$ is any point in $J^*,$ then there is a sequence
$\{y_n\}\subset W^s(y), y_n \rightarrow x.$ In particular, $G^-(y_n)>c$ for all
large $n$, so $G^-(x)\geq c,$ a contradiction since $G^- \equiv 0$ on $J.$

Hence, by Proposition 6, 
there exists a $p \in J^*$ and a nonempty connected component
$V \subset \CC$ of $\CC\setminus J^*_{p,s}$ on which
$G^- \equiv 0.$ Then $V$ is simply connected by Lemma 5: In fact,
if $V$ contains a simple closed curve $\gamma$
whose interior $U$ contains $x\in J^*$, then $G^- \equiv 0$
on $U$ by the maximum principle. Hence by the invariance of $G^-$
it follows that $G^- \equiv 0$ on $H^{-n}(U)$ for all $n \geq 1.$
This contradicts that $H^{-n}(U)$ has to be an unbounded sequence
clustering all over $J^+.$ It follows that there is a weakly stably exposed
point in $J^*.$ Similarly we can find a point $q\in J^*$ which is weakly unstably exposed. \\

\CQFD

\medskip

Next we prove the Theorem 3.

{\bf Proof:}
Pick a weakly stably exposed point $p\in J^*.$ Then in a neighborhood
of $p,$ the unstable lamination contains a laminated hypersurface
with $J^*$ on one side. Let $D$ denote a small unstable disc
centered at $p.$ The forward iterates of this disc become dense in
$J^*.$ Moreover, for every point $q\in H^n(D)\cap J^*$ there is a connected
component $V_{s,q} \subset \CC \setminus J^*_{q,s}$ obtained by following $V_{s,p}$
using the local product structure of $J^*.$ Furthermore, by the local product structure
each such $V_{s,q}$ is simply connected. Hence, arbitrarily close to any point in $J^*$ there is
a laminated hypersurface whis is contained in the local unstable set and
with $J^*$ on one side. 
One can do the same with a weakly unstably exposed point in $J^*$.
Their intersection points give rise to a dense
collection of distinguished boundary points.\\

\CQFD
 
 Next we prove the Main Theorem.
 
 {\bf Proof:}
 We assume that $J \neq J^*$. The hypothesis that $H$ is not volume preserving is
 only used at the end. As in the above proof, for every point $z\in H^n(D) \cap J^*$
 there is a simply connected component of $W^s_z \setminus J^*.$
 Using the local product structure of $J^*$
 for a complex H\'enon map which is hyperbolic on $J^*$, we get that in fact
 for every point $z\in J^*$, there is a simply connected component of
 $W^s_z\setminus J^*.$ Also all connected
 components on which $G^-$ has a zero are simply connected.
 By \cite{BS1998}, if $G^-$ is not identically zero on some
 simply connected component, then $H$ is stably connected. In particular
 all connected components of all $W^s_x\setminus J^*$ for all $x\in J^*$
 are simply connected and all $W^s_x\setminus J^*$ contain at least
 one component where $G^+>0$ and at least one has a component
 where $G^-=0.$ The same applies to the unstable set.

 \begin{lemma}
 Suppose that $H$ is a complex H\'enon map which is hyperbolic on $J^*.$
 Suppose that $J \neq J^*.$
 If $H$ is not stably connected, then each stable manifold $W^s_z\setminus J^*$
 considered as parametrized by $\CC$ 
 contains simply connected components in $\CC$ on which $G^-\equiv 0$ and all
 such components are bounded in $\CC.$
 \end{lemma}
 
 {\bf Proof:}
  Each $W^s_z \setminus J^*$ contains connected components on which $G^->0.$
 However none of these are simply connected, hence they contain curves $\gamma$
 surrounding a nonempty part of $J^*.$ Any simply connected component
 of some $W^s_x\setminus J^*$ can be followed to any other $W^s_z\setminus J^*$ using
 the local product structure and in particular can be followed into the interior of
 such a $\gamma.$ However, unbounded components stay unbounded when
 they are followed in this way. This shows that all simply connected components
 are bounded.\\
 
 \CQFD
 
 \begin{lemma}
 Let $H$ be a complex H\'enon map which is hyperbolic on $J^*$ and suppose that
 $J \neq J^*.$  Then
 $H$ is stably connected and unstably connected.
 \end{lemma}
 
 {\bf Proof:}
 We show that $H$ is stably connected. Unstable connectedness is proved similarly.
 Suppose that $H$ is not stably connected. We pick a point $z\in J^*$ and a
 bounded simply connected component $V$ of $W^s_z\setminus J^*.$
 After changing $z$ if necessary we can pick a curve $\lambda=\phi(t), t\in [0,1]$ with $\lambda(t)\in V, t>0$ and $\lambda(0)=z.$
 We may assume that
 $V\subset W^s_{\epsilon, loc}(J^*)$ for some small $\epsilon>0$ after
 forward iteration if necessary. 
 Next, consider the backward iterates, $V_n:=H^{-n}(V), \gamma_n=H^{-n}(\gamma),
 z_n=H^{-n}(z).$ For $n \geq n_0$ there is a unique point $w_n=H^{-n}(\phi(t_n))\in \gamma_n$
 so that $w_n \in \partial W^s_{\epsilon,loc}(J^*)$ and the curve $\gamma_n$
 is contained in $W^s_{\epsilon,loc}(J^*)$ for $0\leq t <t_n.$ 
 
 \medskip
 
 Let $w_0$ be a cluster point of the sequence $\{w_n\}.$
 Then $w_0$ is in some $W^s_\eta, \eta\in J^*$ and belongs to some
 simply connected component $U$ of $W^s_\eta\setminus J^*$ since $G^-(\eta)=0$
 by continuity. Since this component
 is bounded, we can after replacing with a forward iterate, assume that
 $U$ is contained in the local stable manifold of $\eta.$ Using the local
 product structure it follows that the diameter
 of $V$ is arbitrarily small, a contradiction.\\
 
 \CQFD
 
 Using \cite{BS1998} Theorem 5.1 and Corollary 7.4 we get:
 
 \begin{theorem}
 Suppose that $H$ is a complex H\'enon map which is hyperbolic on $J^*.$
 Moreover assume that $J\setminus J^* \neq \emptyset.$
 Then $H$ is stably and unstably connected. Morever $J$ is connected
 and $H$ is volume preserving. 
 \end{theorem}
 
 The Main Theorem is now an immediate consequence.\\
 
 \CQFD

\end{document}